# Fourier series and periodicity

Donal F. Connon

8 December 2014

**Abstract**

A large number of the classical texts dealing with Fourier series more or less state that the hypothesis of periodicity is required for pointwise convergence. In this paper, we highlight the fact that this condition is not necessary.

**1. Whither periodicity?**

When I first studied Fourier series some 40 years ago, as part of an applied mathematics course, I was told that it was <u>necessary</u> for the function to be periodic, i.e. $f(x) = f(x + p)$ where (usually) $p = 2\pi$. If $f(x)$ was defined on an interval, say, $[-\pi, \pi]$, we were told, as in [13, p.244] and [12, p.299], to extend $f(x)$ periodically (and also to graphically display the extended function so as to aid visualisation).

One would then obtain a Fourier series

$$f(x) = \frac{1}{2}a_0 + \sum_{n=1}^{\infty}(a_n \cos nx + b_n \sin nx)$$

and, in view of the fact that both $\cos nx$ and $\sin nx$ were $2\pi$ periodic, it is therefore clear that we do indeed end up with a representation of $f(x)$ that is periodic.

However, as I mentioned in [7], I was never able to truly fathom why such periodicity should, ab initio, be a requirement in the Dirichlet conditions used in the rigorous treatment of Fourier series (see, for example, Apostol [1, p.319], Bressoud [3, p.223], Titchmarsh [14, p.401] and Whittaker and Watson [16, p.164]). Why indeed should the nature of a real function outside of its defined finite domain have any effect on the pointwise convergence of the corresponding Fourier series?

Apostol [1, p.323] stated that "The hypothesis of periodicity, which appears in all the convergence theorems dealing with Fourier series, is not as serious a restriction as it may appear to be at first sight. If a function is initially defined on a finite interval, say $[a,b]$, we can always extend the definition of $f$ outside $[a,b]$ by imposing some sort of periodicity condition. For example, if $f(a) = f(b)$, we can define $f$ everywhere on $(-\infty, +\infty)$ by requiring the equation $f(x) = f(x + p)$ to hold for every $x$ where $p = b - a$. (The condition $f(a) = f(b)$ can always be brought about by changing the value of $f$ at one of the endpoints if necessary. This does not affect the existence or the values of the integrals which are used to compute the Fourier coefficients of $f$.)"

What Apostol says is of course perfectly correct, but we are still in the dark as to why it is considered that periodicity is required ab initio. We do not really fare any better by examining other classical texts such as Zygmund [17, p.8] where it is simply stated that "…, it is *convenient* to assume (as we shall always do in what follows) that the functions whose Fourier series we consider are defined not only in an interval of length $2\pi$ but for all values of $x$, by the condition of periodicity $f(x+2\pi) = f(x)$." [emphasis added]. Hobson [10, p.490] makes a similar comment regarding periodicity and convenience.

Tolstov [15] starts his book with a section on periodic functions and appears to require the periodicity condition throughout his text. This is also the approach taken, inter alia, by Bachman, Narici and Beckenstein [2, p.184], Churchill [6, p.90], Knopp [11, p.358] and Whittaker and Watson [16, p.164].

**2. Demonstration that periodicity is not required**

The most fundamental formula in the theory of Fourier series is Dirichlet's integral which may be stated [5, p.227] for $0 < a < \pi$

(1) $$\lim_{\mu \to \infty} \int_0^a f(x) \frac{\sin \mu x}{\sin x} dx = \frac{\pi}{2} f(0+)$$

where $\mu$ is not necessarily an integer.

Let us first examine how Apostol [1, p.317] approaches the theory of Fourier series where he considers the series expansion on the interval $[0, 2\pi]$. We have the partial sum

$$s_n(x) = \frac{1}{2}a_0 + \sum_{k=1}^n (a_k \cos kx + b_k \sin kx)$$

where the Fourier coefficients are

$$a_k = \frac{1}{\pi} \int_0^{2\pi} f(t) \cos kt \, dt \qquad b_k = \frac{1}{\pi} \int_0^{2\pi} f(t) \sin kt \, dt$$

Hence we have

$$s_n(x) = \frac{1}{\pi} \int_0^{2\pi} f(t) \left[ \frac{1}{2} + \sum_{k=1}^n (\cos kt \cos kx + \sin kt \sin kx) \right] dt$$

$$= \frac{1}{\pi} \int_0^{2\pi} f(t) \left[ \frac{1}{2} + \sum_{k=1}^n \cos k(t-x) \right] dt$$

$$= \frac{1}{\pi} \int_0^{2\pi} f(t) \frac{\sin\left[(2n+1)(t-x)/2\right]}{2\sin\left[(t-x)/2\right]} dt$$

so that

(4) $$s_n(x) = \frac{1}{\pi} \int_0^{2\pi} f(t) D_n(t-x) \, dt$$

where the function $D_n$ is called Dirichlet's kernel.

Since $D_n$ is periodic with period $2\pi$, and postulating that $f$ is also periodic with period $2\pi$, we can replace the interval of integration by $[x-\pi, x+\pi]$ and then make a translation $u = t - x$ to get

$$s_n(x) = \frac{1}{\pi} \int_{x-\pi}^{x+\pi} f(t) D_n(t-x) \, dt$$

$$= \frac{1}{\pi} \int_{-\pi}^{\pi} f(x+u) D_n(u) \, du$$

Since $D_n(u) = D_n(-u)$ we obtain (on the assumption that $f$ is periodic)

(5) $$s_n(x) = \frac{2}{\pi} \int_0^{\pi} \frac{f(x+t) + f(x-t)}{2} D_n(t) \, dt$$

which we may write as

$$= \frac{4}{\pi} \int_0^{\pi/2} \frac{f(x+2v) + f(x-2v)}{2} D_n(2v) \, dv$$

Then applying (1) we obtain for $0 < x < 2\pi$

$$\lim_{n \to \infty} s_n(x) = \frac{1}{2} [f(x-0) + f(x+0)]$$

Enlightenment of the fact that the periodicity condition is redundant finally came through reading Carslaw's excellent treatise, Introduction to the theory of Fourier's series and integrals [5, p.230].

Instead of assuming that $f$ is periodic, we start with (4)

$$s_n(x) = \frac{1}{\pi} \int_0^{2\pi} f(t) D_n(t-x) \, dt$$

$$= \frac{1}{\pi} \int_0^x f(t) D_n(t-x) \, dt + \frac{1}{\pi} \int_x^{2\pi} f(t) D_n(t-x) \, dt$$

and make the substitution $t = x - 2v$ in the first integral and $t = x + 2v$ in the second integral to obtain

(6) $$s_n(x) = \frac{2}{\pi} \int_0^{\frac{x}{2}} f(x-2v) D_n(2v)\, dv + \frac{2}{\pi} \int_0^{\pi - \frac{x}{2}} f(x+2v) D_n(2v)\, dv$$

Then using (1) we again obtain for $0 < x < 2\pi$

(7) $$\lim_{n \to \infty} s_n(x) = \frac{1}{2}[f(x-0) + f(x+0)]$$

It is thus seen that there is no need to assume ab initio that $f$ is periodic.

When $x = 0$ we have

$$s_n(0) = \frac{2}{\pi} \int_0^{\pi} f(2v) D_n(2v)\, dv$$

and from (8) below we find that

$$\lim_{n \to \infty} s_n(0) = \frac{1}{2}[f(0+) + f(2\pi-)]$$

Similarly, when $x = 2\pi$ we have

$$s_n(2\pi) = \frac{2}{\pi} \int_0^{\pi} f(2\pi - 2v) D_n(2v)\, dv$$

and (8) tells us that

$$\lim_{n \to \infty} s_n(0) = \frac{1}{2}[f(0+) + f(2\pi-)]$$

### 3. Extension of Dirichlet's integral

As noted by Bromwich [4, p.492] we may increase the range of Dirichlet's integral (1) to $\pi$ in the case where $\mu = 2N + 1$ and $N$ is an integer. We have

$$\int_0^{\pi} f(x) \frac{\sin(2N+1)x}{\sin x}\, dx = \int_0^{\pi/2} f(x) \frac{\sin(2N+1)x}{\sin x}\, dx + \int_{\pi/2}^{\pi} f(x) \frac{\sin(2N+1)x}{\sin x}\, dx$$

and make the substitution $x = \pi - t$ in the second part. This immediately gives us

$$\int_0^\pi f(x)\frac{\sin(2N+1)x}{\sin x}dx = \int_0^{\pi/2}[f(x)+f(\pi-x)]\frac{\sin(2N+1)x}{\sin x}dx$$

and employing (1) we obtain

(8) $$\lim_{N\to\infty}\int_0^\pi f(x)\frac{\sin(2N+1)x}{\sin x}dx = \frac{\pi}{2}[f(0+)+f(\pi-)]$$

In fact, as noted below, we may extend the range of the integration well beyond $\pi$. For example, we consider

$$\int_0^{2\pi} f(x)\frac{\sin(2N+1)x}{\sin x}dx = \int_0^\pi f(x)\frac{\sin(2N+1)x}{\sin x}dx + \int_\pi^{2\pi} f(x)\frac{\sin(2N+1)x}{\sin x}dx$$

and we write

$$\int_\pi^{2\pi} f(x)\frac{\sin(2N+1)x}{\sin x}dx = \int_\pi^{3\pi/2} f(x)\frac{\sin(2N+1)x}{\sin x}dx + \int_{3\pi/2}^{2\pi} f(x)\frac{\sin(2N+1)x}{\sin x}dx$$

The substitution $x = \pi + t$ gives us

$$\int_\pi^{3\pi/2} f(x)\frac{\sin(2N+1)x}{\sin x}dx = \int_0^{\pi/2} f(\pi+t)\frac{\sin(2N+1)t}{\sin t}dt$$

and with $x = 2\pi - t$ we obtain

$$\int_{3\pi/2}^{2\pi} f(x)\frac{\sin(2N+1)x}{\sin x}dx = \int_0^{\pi/2} f(2\pi-t)\frac{\sin(2N+1)t}{\sin t}dt$$

Therefore, using Dirichlet's integral (1) we see that

$$\lim_{N\to\infty}\int_0^{2\pi} f(x)\frac{\sin(2N+1)x}{\sin x}dx = \frac{1}{2}[f(0+)+f(\pi-)+f(\pi+)+f(2\pi-)]$$

More generally we see that

$$\lim_{N\to\infty}\int_0^{m\pi} f(x)\frac{\sin(2N+1)x}{\sin x}dx = \frac{1}{2}[f(0+)+f(m\pi-)] + \frac{1}{2}\sum_{n=1}^{m-1}[f(n\pi-)+f(n\pi+)]$$

and, if $f(x)$ is continuous, this may be written as

$$\lim_{N\to\infty} \int_0^m f(x) \frac{\sin(2N+1)\pi x}{\sin \pi x} dx = \frac{1}{2}[f(0)+f(m)] + \sum_{n=1}^{m-1} f(n)$$

In [8] we used this to obtain a version of the Poisson summation formula

$$\sum_{n=0}^m f(n) = \int_0^m f(x)\,dx + \frac{1}{2}[f(0)+f(m)] + 2\sum_{n=1}^\infty \int_0^m f(x)\cos 2n\pi x\,dx$$

Letting $m \to \infty$, and assuming that $\lim_{n\to\infty} f(n) = 0$ (which is of course required for the convergence of the infinite series), we obtain

$$\frac{1}{2}f(0) + \sum_{n=1}^\infty f(n) = \int_0^\infty f(x)\,dx + 2\sum_{n=1}^\infty \int_0^\infty f(x)\cos 2\pi nx\,dx$$

We may also note that

$$\int_0^a f(x) \frac{\sin(2N+1)x}{\sin x} dx = \int_0^a f(x) \frac{\sin 2Nx \cos x}{\sin x} dx + \int_0^a f(x)\cos 2Nx\,dx$$

and, for suitable functions, the Riemann Lebesgue lemma tells us that

$$\lim_{N\to\infty} \int_0^a f(x)\cos 2Nx\,dx = 0$$

and thus for $0 < a < \pi$

$$\lim_{N\to\infty} \int_0^a f(x)\cot x \sin 2Nx\,dx = \frac{\pi}{2} f(0+)$$

**4. Open access to our own work**

This paper contains references to various other papers and, rather surprisingly, all of them are currently freely available on the internet. Surely now is the time that <u>all</u> of <u>our</u> work should be freely accessible by <u>all</u>. The mathematics community should lead the way on this by publishing <u>everything</u> on arXiv, or in an equivalent open access repository. We think it, we write it, so why hide it? You know it makes sense.

Wessex House,
Devizes Road,
Upavon,
Pewsey,
Wiltshire SN9 6DL
dconnon@btopenworld.com